\documentclass[fleqn]{article}
\usepackage{latexsym}
\usepackage{amstex}

\def\BibTeX{{\rm B\kern-.05em{\sc i\kern-.025em b}\kern-.08em T\kern-.1667em\lower.7ex\hbox{E}\kern-.125emX}}

\begin{document}

\title{A Simple Proof of the Four-Color Theorem}
\author{Fayez A. Alhargan\footnote{KACST, CERI, Riyadh, Saudi Arabia, Email: alhargan at   kacst.edu.sa}}

\maketitle

\begin{abstract}
A simpler proof of  the four color theorem is presented. The proof was reached using a series of equivalent theorems. First the maximum number of edges of a planar graph is obatined as well as the minimum number of edges for  a complete graph. Then it is shown that for the theorem to be false there must exist a complete planar graph of $h$ edges such that $h>4$. Finally the theorem  is proved to be true by showing that there does not exist a complete planar graph with $h>4$. 
\end{abstract}

\section{Introduction}
The Four Color Problem has been investigated by many mathematicians since 1850 with no conclusive mathematical proof \cite{Saaty} and \cite{Fritsch}. A proof was demonstrated by Kenneth Appel and Wolfgang Haken in 1976 using very sophisticated computer program \cite{Appel1},\cite{Appel2} and \cite{Appel3}. The the Appel-Haken proof has not been fully accepted \cite{Robertson}, as one part of the proof relies on a computer which cannot be verified by hand and the other part of the proof is very complicated and tedious to verify by hand. Robertson et. al. \cite{Robertson} have proposed a simpler proof using similar methodology to that of Appel-Haken, which still required a computer, but the part which required hand verification was not as tedious.  In this paper the truth of the theorem is shown by proving  that there does not exist a complete planar subgraph $H=S(G)$, which consist of more than four vertices. For accessibility to wider audience the proofs are presented in away that can be accessible to mathematics undergraduates.

\section{Statement of Main Results}
Definitions:\newline
$M$: is a general map.\newline
$G$: is the underlying graph of $M$,  $G=U(M)$.\newline
$H$: is a subgraph of $G$ with $h$ vertices, $H=S(G)$.\newline
A planar graph: is a graph which can be drawn in the plane so that its edges do not cross.\newline
A complete graph: is a graph which has its vertices  connected in such a way  that every vertex is connected to every other vertex.\newline

A complete graph is shown in figure 1. However, the graph is not planar as will be shown later, whereas figure 2 is planar but not complete.\newline

\begin{figure}
\begin{center}
\begin{picture}(200,200)(-40,0)
\put(60,60){\circle*{4}}
\put(55,65){1}
\put(110,110){\circle*{4}}
\put(105,115){2}
\qbezier(60,60)(75,75)(110,110)
\put(10,10){\circle*{4}}
\put(0,0){3}
\qbezier(10,10)(60,60)(110,110)
\put(10,110){\circle*{4}}
\put(5,115){4}
\qbezier(10,110)(60,110)(110,110)
\qbezier(10,110)(10,60)(10,10)
\qbezier(10,110)(35,75)(60,60)
\qbezier(110,110)(75,35)(10,10)
\put(110,10){\circle*{4}}
\put(105,0){5}
\qbezier(110,10)(60,10)(10,10)
\qbezier(110,10)(110,60)(110,110)
\qbezier(110,10)(200,200)(10,110)
\qbezier(110,10)(75,35)(60,60)
\end{picture}
\end{center}
\label{complgrf}
\caption{Complete graph}
\end{figure}
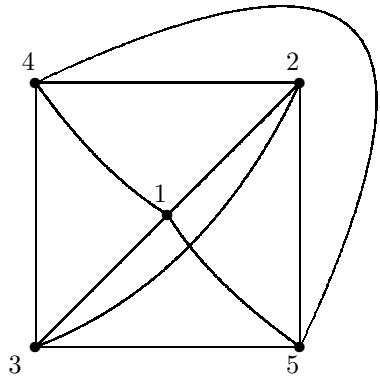

\begin{figure}
\begin{center}
\begin{picture}(200,200)(-40,0)
\put(60,60){\circle*{4}}
\put(55,65){1}
\put(110,110){\circle*{4}}
\put(105,115){2}
\qbezier(60,60)(75,75)(110,110)
\put(10,10){\circle*{4}}
\put(0,0){3}
\qbezier(10,10)(60,60)(110,110)
\put(10,110){\circle*{4}}
\put(5,115){4}
\qbezier(10,110)(60,110)(110,110)
\qbezier(10,110)(10,60)(10,10)
\qbezier(10,110)(35,75)(60,60)
\qbezier(110,110)(75,35)(10,10)
\put(110,10){\circle*{4}}
\put(105,0){5}
\qbezier(110,10)(60,10)(10,10)
\qbezier(110,10)(110,60)(110,110)
\qbezier(110,10)(200,200)(10,110)
\end{picture}
\end{center}
\label{plangrf} 
\caption{Planar graph}
\end{figure}
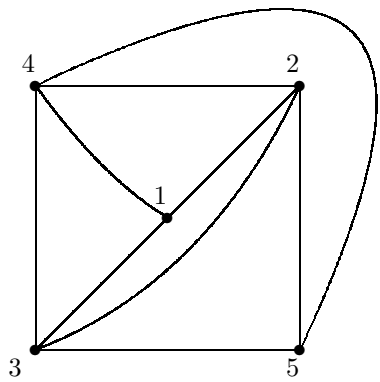

Postulate: Maps consisting of regions and boundaries can be represented by a dual graph consisting of vertices and edges. Where vertices represent regions and edges represent boundaries.

Theorem $C1$: Four colors are sufficient to color any map drawn in a plane so that no two regions with a common boundary line are colored with the same color.

Theorem $C2$: Four colors are sufficient to color all vertices of a planar graph such that no two vertices connected by an edge are colored with the same color.

Theorem $C3$: There does not exist a subgraph  $H=S(G)$, with $h$ vertices; such that $H$ is a complete planar graph and $h>4$.

Theorem T1: $ C2 \Leftrightarrow  C1 $\newline
Proof of T1: {see (\cite{Saaty},p.16)}

Theorem T2: $ C3 \Leftrightarrow C2$ \newline
Proof of T2:\newline
If $C2$ is false, this implies that there exist a complete planar subgraph  $H=S(G)$, such  that $H$ has $h$ vertices and $h>4$. Only the existence of such a graph would make $C2$ false. Moreover, if such a graph exist, then $C3$ is false. \newline
Hence, $ \sim C3 \Leftrightarrow \sim C2$. \newline
Therefore, $ C3 \Leftrightarrow C2 $

Theorems T1 and T2 imply that the three color theorems $C1$, $C2$ and $C3$ are equivalent theorems, proving any one of them proofs them all.

Theorem T3:Let $M$ be a map whose underlying graph $G$ has $n$ vertices and $m$ edges and suppose $M$  has $r$ regions. Then $n-m+r=2$.\newline
Proof of T3: {see (\cite{Saaty},p.18)}

Corollary 1:\newline
Let $G$ be a simple connected planar graph with $n$ vertices and $m$ edges. Then
\begin{equation}
  \label{c1}
 m \leq 3n-6, \quad n\ge 3  
\end{equation}

Proof of Corollary 1 (\cite{Saaty},p.19) : \newline
Let $M$ be a map with $r$ regions such that $U(M)=G$. Since $G$ is simple, every region of $M$ has at least three sides. Moreover, every edge of $M$ lies on the boundary of at most two regions. \newline
Therefore, $ 3r \leq 2m $. \newline
But from theorem T3, $ r=2-n+m $.   Substituting for $r$, gives,  $ m \leq 3n -6 $

Corollary 2:\newline
If $G$ is a complete graph with $n$ vertices and $m$ edges. Then 
\begin{equation}
  \label{c2}
  m = \sum^{n-1}_{r=0} r = \frac{1}{2} (n-1) n 
\end{equation}

Proof of Corollary 2:\newline
Starting with a single vertex, then $n=1$ and $m=0$. Adding to the graph a second vertex and connecting it to all the vertices of the graph, gives $n=2$ and  $m=0+1$, then  adding a third vertex, gives $n=3$ and $m=0+1+2$, also adding a fourth vertex, gives $n=4$ and $m=0+1+2+3$. Repeating this process $n$ times will result in a graph which is complete (fully connected), having $n$ vertices and $m=0+1+2+3+...+(n-1)$ edges.

Proof of the color theorem $C3$:\newline
Equating (\ref{c1}) and (\ref{c2}); gives,
\begin{equation}
 \frac{1}{2} (n-1) n \leq 3(n-2) 
\end{equation}
i.e.
\begin{equation}
 n^2 - 7n +12 \leq 0 
\end{equation}
or
\begin{equation}
  \label{c3e2}
  (n-3) (n-4) \leq 0 
\end{equation}
Solving equation (\ref{c3e2}), gives
\begin{equation}
  \label{c3e3}
 3 \leq n \leq 4 
\end{equation}
From (\ref{c3e3}), the maximum possible number of vertices for a complete planar graph is 4. 

Now, from  equation (\ref{c3e3}); for a subgraph $H$ of $h$ vertices with $h>4$, there  does not exist $H=S(G)$, such that  $H$ is a complete  planar graph. \newline
Hence C3 is true. \newline
But, $ C3 \Leftrightarrow C2$ and  $ C2 \Leftrightarrow C1 $, 
Since $C3$ is true, then by theorem T2, $C2$ is true, and in turn by theorem T1, $C1$ is true. Hence the four color theorem is true.

Finally, can we conjure  about the minimum number of colors required for the other dimensions. For instance the zero dimension (i.e. a point), one color is sufficient, for one dimension (i.e.  line segments) two colors are sufficient and for two dimensions (i.e. planar maps) as has been proved, four colors are sufficient. This suggests that for $m$-dimensions, the minimum number of colors required are $m=2^n$. Is this expression true for three dimensions and higher dimensions?.

\section{Conclusion}
The proof can be readily deduced if one considers what type of a map is needed to make the theorem false. Thinking along this line leads to a map which consist of at least five regions and to make the theorem false, the regions must touch each other completely (i.e. every region must touch every other region). However trying to realize such a map, it becomes apparent  that such a map can not be realized on the plane. Hence  the theorem must be true.

\providecommand{\bysame}{\leavevmode\hbox to3em{\hrulefill}\thinspace}
\providecommand{\MR}{\relax\ifhmode\unskip\space\fi MR }
\providecommand{\MRhref}[2]{%
  \href{http://www.ams.org/mathscinet-getitem?mr=#1}{#2}
}
\providecommand{\href}[2]{#2}

\end{document}